\documentclass[11pt,twoside]{amsart} 

\title{Equivariant stratifold homology theories}
\author{Julia Weber}


\usepackage[latin1]{inputenc} %% deutsche Umlaute normal tippen
\usepackage{bbm}
\usepackage{amssymb}
\usepackage{amsthm}
\usepackage{xspace}
\usepackage[all]{xy}\CompileMatrices

\newcommand{\Z}{\ensuremath{\mathbb{Z}}\xspace}
\newcommand{\R}{\ensuremath{\mathbb{R}}\xspace}
\newcommand{\Nat}{\ensuremath{\mathbb{N}}\xspace}
\newcommand{\Zz}{\ensuremath{\Z/2}\xspace}
\newcommand{\N}{\ensuremath{\mathcal{N}}\xspace}

\newcommand{\Cl}{\ensuremath{\mathcal{C}}\xspace}

\newcommand{\F}{\ensuremath{\mathcal{F}}\xspace}
\renewcommand{\H}{\ensuremath{\mathcal{H}}\xspace}

\newcommand{\oS}{\ensuremath{\text{\emph{\r{S}}}}\xspace}

\newcommand{\isopfeil}{\ensuremath{\xrightarrow{\sim}}\xspace}

\newtheorem{Def}{Definition}[section]

\newtheorem{Cor}[Def]{Corollary}
\newtheorem{Prop}[Def]{Proposition}
\newtheorem{Lem}[Def]{Lemma}

\newcommand{\Rem}{\noindent \bfseries Remark: \mdseries}
\newcommand{\Rems}{\noindent \bfseries Remarks: \mdseries}
\begin{document}

\maketitle

\begin{abstract}

We define equivariant homology theories using bordism of stratifolds with a $G$-action, where $G$ is a discrete group. Stratifolds are a generalization of smooth manifolds which were introduced by Kreck~\cite{kreck}. He defines homology theories using bordism of suitable stratifolds. We develop the equivariant generalization of these ideas. 

\end{abstract}

\section{Introduction and definitions} \label{sec1}

Let $G$ be a discrete group. In this paper, we are going to define equivariant homology theories using bordism of stratifolds with a $G$-action. Stratifolds are defined inductively. An $n$-dimensional stratifold is obtained by gluing an $n$-dimensional smooth manifold to an $n-1$-dimensional stratifold via its boundary. The manifolds are equipped with collars which yield a smooth structure on the stratifolds. 

Homology theories can be defined via bordism of certain classes of stratifolds. In particular, one obtains singular homology as a bordism theory. Stratifolds and the corresponding homology theories were developed by Kreck~\cite{kreck}, where one can find details concerning the relevant constructions. 

We are concerned with the development of corresponding equivariant homology theories, a question raised by Kreck. We define the notion of $G$-stratifolds, for discrete groups $G$, and prove that they have the technical properties necessary to define a $G$-homology theory. We also treat an induction structure, thereby obtaining an equivariant homology theory. We give several examples of equivariant homology theories thus obtained. In particular, we have a generalization of singular homology, an equivariant theory with coefficients in the Burnside ring.

%The reader is assumed to know about the definition of stratifolds and to be familiar with the construction of homology theories using bordism of stratifolds. 

\begin{Def}
A stratifold $S$ is called a \emph{$G$-stratifold} if there is a smooth $G$-action on $S$. This means there is a smooth map of stratifolds (the discrete group $G$ is a $0$-dimensional smooth manifold, so $G\times S$ with the product stratifold structure is a stratifold)
\begin{eqnarray*}
\theta : G\times S & \to & S\\
 (g,s) & \mapsto & gs
\end{eqnarray*}
such that
\begin{enumerate}
\item es=s
\item (gh)s=g(hs).
\end{enumerate}
$S$ is called a \emph{proper} $G$-stratifold if the action of $G$ on $S$ is proper, which means that the map $G\times S\to S\times S, (g,s)\mapsto (gs,s)$ is proper.
\end{Def}

\Rem 
The map $\theta_g:S\to S, s\mapsto gs$ is in $Aut(S)$ for all $g\in G$.

This means that $g\oS_i=\oS_i$ for all $g\in G$, that collars are preserved, collars of the strata as well as collars of $S$, and that the collar maps are equivariant (with respect to the trivial $G$-action on the interval). Thus we can also define the operation of $G$ on each stratum seperately (and not only on its interior $\oS_i \hookrightarrow S$) by continuing it onto the boundary.
This leads to the observation that we could also have defined a $G$-stratifold by saying: The strata $S_i$ are $G$-manifolds with $G$-equivariant collars which are glued together with $G$-equivariant maps.

\begin{Def}
Let $X$ be a $G$-space. An $n$-dimensional \emph{singular $G$-stratifold} is a pair $(S,f)$, where $S$ is an $n$-dimensional cocompact $G$-stratifold (i.e. $G\backslash S$ is compact) and $f:S\to X$ is a $G$-equivariant map. 
\end{Def}

As for usual homology theories defined via stratifolds, we can restrict ourselves to stratifolds satisfying certain conditions, thereby obtaining many different $G$-homology theories depending on the conditions imposed. We call a class of stratifolds a bordism class if it has properties allowing an equivalence relation ``bordism'' of its stratifolds and if the restrictions are sufficient to ensure that this bordism theory is a homology theory.

\begin{Def}\label{defbordclass}
A class $\Cl$ of stratifolds is called a \emph{bordism class} if it satisfies the following axioms:

{\bf Axiom 1:} If $S$ is a closed $G$-stratifold in $\Cl$, $S\times I$ (where the smooth manifold $I=[0,1]$ with boundary has trivial $G$-operation) is in $\Cl$.

{\bf Axiom 2:} If $S$ and $S'$ are cocompact $G$-stratifolds in $\Cl$, then $S\amalg S'$ is in $\Cl$. If $T'$ and $T''$ are cocompact $G$-stratifolds with boundary such that $\partial T'=S\cup S'$ with $S$ and $S'$ $G$-invariant and $\partial T''=S'\cup S''$ with $S'$ and $S''$ $G$-invariant and such that $S\cap S'=S'\cap S''=\partial S=\partial S'=\partial S''$ is an $(n-2)$-dimensional stratifold then $T'\cup_{S'}T''\in \Cl$, and $\partial(T'\cup_{S'}T'')=S\cup_{\partial S'}S''\in \Cl$. 

{\bf Axiom 3:} If $S$ is a cocompact $n$-dimensional $G$-stratifold in $\Cl$, $\rho:S\to \R$ is a smooth $G$-invariant map and $t\in \R$ is a regular value of $\rho$, then $\rho^{-1}(t)$ is in $\Cl$, and $\rho^{-1}([t,\infty))$ and $\rho^{-1}((-\infty,t])$ are in $\Cl$. 
\end{Def}

\Rems
\begin{itemize}
\item Except for the fact that the boundary components in Axiom 2 have to be $G$-invariant and that $\rho:S\to \R$ in Axiom 3 is $G$-invariant, the $G$-operation does not play a role in the formulation of these axioms, they are completely analogous to the ones used for usual homology theories.
\item It is shown in~\cite{kreck} that in the situation of Axiom 3 the pre-image of a regular value, $\rho^{-1}(t)$, is an $(n-1)$-dimensional stratifold, so to check Axiom 3 it is only necessary to verify that this fulfills the restrictions again.
\item It is useful to observe that these axioms are usually fulfilled for conditions given ``locally'', that is in terms of a neighborhood of $x$ for all $x\in S$. In particular, regular ($G$-)stratifolds are defined by the fact that each $x\in S$ has a neighborhood $V_{U_x}\cong B^i\times F$. A $G$-invariant neighborhood is then given by 
$$G V_{U_x} \cong G(B^i\times F)=(GB^i)\times F=(\cup_{g\in G/G_x}gB^i)\times F,$$ 
since we have equivariant collars. (Remember that $F\cong r_i^{-1}(\{x\})$, where $r_i$ is the projection along the collars.) If we only impose conditions on $F$, this gives a class of stratifolds satisfying all the axioms. (An example for this are Euler stratifolds. In their definition we require $F$ to fulfill the condition $\chi(F\backslash F^{(0)})\equiv 0 \mod 2 $.)
\item Another case where the axioms are automatically satisfied is if we impose conditions on the emptyness/non-emptyness of strata of certain codimensions - the most prominent example is the special case ``codimension 1-stratum empty'' which leads to singular homology with $\Zz$-coefficients~\cite{kreck}.
\item If we want to impose orientability conditions (as for example the extra condition ``top-dimensional stratum oriented'' leading to singular homology with $\Z$-coefficients instead of $\Zz$-coefficients), we have to require all $G$-operations to be orientation preserving. 
\end{itemize}

\begin{Def}
Two $n$-dimensional cocompact singular $G$-stratifolds $(S,f)$ and $(S',f')$ in $\Cl$ are called \emph{bordant} (in $\Cl$) $\Leftrightarrow$ there is an $(n+1)$-dimensional cocompact $G$-stratifold $T$ in $\Cl$ with boundary and an equivariant map $F:T\to X$ such that $\partial T=S\amalg S', F|_{\partial T}=f\amalg f'$.
\end{Def}

\begin{Prop}
The relation ``bordant'' (in $\Cl$) defined above is an equivalence relation.
\end{Prop}
\begin{proof}
Reflexivity: Clear. (Take $T:=S\times [0,1]$, with trivial $G$-operation on $[0,1]$.) \\
Symmetry: Clear.\\
Transitivity: If $(S,f) \sim_{(T,F)} (S',f')$ and $(S',f') \sim_{(T',F')} (S'',f'')$, then we can glue together $T$ and $T'$ along their common boundary component $S'$ to obtain a bordism between $(S,f)$ and $(S'',f'')$ given by $(T\cup_{S'}T', F\cup_{f'}F')$. (Remember that we have an equivariant collar.)
\end{proof}

\Rem
Under disjoint union the equivalence classes $[S,f]$ with respect to the equivalence relation ``bordant'' (in $\Cl$) of $n$-dimensional cocompact singular $G$-stratifolds $(S,f)$ (in $\Cl$) form a group; the zero element is $[\emptyset]$, each element is its own inverse. Thus we get a $\Zz$-vector space. (The question whether this is a set has to be resolved analogously to usual stratifold bordism~\cite{kreck}.)\\

We want to show that these bordism groups form a $G$-homology theory. We first need some technical statements:

\begin{Lem} \label{partitionofunity}
Let $M$ be a proper $G$-manifold. Then for every $G$-invariant open covering $\{V_\alpha\}_{\alpha\in A}$ there is a subordinate smooth $G$-invariant partition of unity. More precisely, there is a collection of smooth $G$-invariant maps $\{\lambda_\beta:M\to (-\varepsilon,1+\varepsilon)\}$ such that
%KANN ICH DAS NOCH AENDERN ZU $[0,1]$?
\begin{enumerate}
\item There is a locally finite open refinement $\{U_\beta\}_{\beta \in B}$ of $\{V_\alpha\}_{\alpha\in A}$ such that $supp(\lambda_\beta)\subseteq U_\beta$ for all $\beta \in B$.
\item $\sum_{\beta\in B} \lambda_\beta(x)=1$ for all $x\in M$.
\end{enumerate}
\end{Lem}
\begin{proof}
$G$ acts properly on $M$, so $G\backslash M$ is a paracompact Hausdorff space. Also we know: For every $x\in M$ there is a smooth slice $S_x$ in $x$~\cite{palais}.
%Illman 1.6, S.137 Palais P3
For every $x\in M$ choose an open $G$-invariant neighborhood $W_x\subseteq G\times_{G_x} S_x$. Let $\pi:M\to G\backslash M$ be the projection. Since $G\backslash M$ is paracompact, there is a locally finite common refinement $\{U'_\beta\}_{\beta \in B}$ of $\{\pi(W_x)\}_{x \in M}$ and $\{\pi(V_\alpha)\}_{\alpha\in A}$ and a partition of unity $\{\lambda'_\beta\}$ subordinate to $\{U'_\beta\}_{\beta \in B}$. Set $U_\beta:=\pi^{-1}(U'_\beta)$ for all $\beta \in B$. This is a locally finite common refinement of $\{W_x\}_{x \in M}$ and $\{V_\alpha\}_{\alpha\in A}$. Define $\lambda''_\beta$ on all of $M$ by setting $\lambda''_\beta(x):=\lambda_\beta'(G x)$ for all $x\in M$. This is a $G$-invariant partition of unity on $M$ subordinate to $\{U_\beta\}_{\beta \in B}$. Now change every $\lambda''_\beta$ into a smooth map: We know that $supp(\lambda''_\beta)\subseteq U_\beta \subseteq W_x\subseteq G\times_{G_x} S_x$ for an $x\in M$. Because of the $G$-invariance $\lambda''_\beta$ is completely defined by its values on $S_x$. Now we approximate $\lambda''_\beta|_{S_x}:S_x\to [0,1]$ by a smooth map $\tilde{\lambda}''_\beta:S_x\to (-\varepsilon,1+\varepsilon)$. We make this map $G$-invariant by setting 
$$\tilde{\lambda}_\beta(s):=\frac{1}{|G_x|} \sum_{g\in G_x} \tilde{\lambda}''_\beta(gs).$$  
%(For $G$ discrete, otherwise use the Haar-Integral on $G_x$; $G_x$ is compact since the operation is proper!) ZITAT!!!\\
This is well defined because $G_x$ is finite since $G$ acts properly on $M$. So now we have a smooth $G_x$-invariant map $\tilde{\lambda}_\beta:S_x\to  (-\varepsilon,1+\varepsilon)$. We define $\tilde{\lambda}''_\beta([g,s]):=\tilde{\lambda}''_\beta(s)$ for all $[g,s]\in G\times_{G_x} S_x$ and continue this map (which has $supp(\tilde{\lambda}''_\beta)\subseteq U_\beta \subseteq W_x$) by $0$ to all of $M$. This is a well defined $G$-invariant map $\tilde{\lambda}_\beta:M\to  (-\varepsilon,1+\varepsilon)$. We do this for all $\beta \in B$ and set $\lambda_\beta(x):=\frac{\tilde{\lambda}_\beta(x)}{\sum_{\beta\in B} \tilde{\lambda}_\beta(x)}.$(The sum in the denominator is finite since the partition of unity $\{U_\beta\}_{\beta \in B}$ was locally finite.)

Then we have $\sum_{\beta\in B} \lambda_\beta(x)=1$ for all $x\in M$, and $\{\lambda_\beta\}$ is the desired smooth $G$-invariant partition of unity. 
\end{proof}

\begin{Lem}
Let $S$ be a proper cocompact $G$-stratifold. Then every $G$-invariant open covering $\{V_\alpha\}_{\alpha\in A}$ has a subordinate smooth $G$-invariant partition of unity.
\end{Lem}
\begin{proof} 
% eigentlich muessten die Argumente mit \beta statt mit \alpha laufen, ist aber noch komplizierter zu formulieren, muesste aber genauso durchgehen
By Lemma~\ref{partitionofunity} every stratum $S_i$ of $S$ has a smooth $G$-invariant partition of unity subordinate to $\{V_\alpha\cap S_i\}_{\alpha\in A}$. With the $G$-equivariant collars one inductively combines these to finally obtain a partition of unity on $S$. The partition of unity on $\{V_\alpha\cap S^{(0)}\}_{\alpha\in A}$ is obvious. So now suppose a smooth $G$-invariant partition of unity $\{(\lambda_{S^{(n-1)}})_\alpha\}$ on $S^{(n-1)}$ subordinate to $\{V_\alpha\cap S^{(n-1)}\}_{\alpha\in A}$ and a smooth $G$-invariant partition of unity $\{(\lambda_{S_n})_\alpha\}$ on $S_n$ subordinate to $\{V_\alpha\cap S_n\}_{\alpha\in A}$ are given. Choose a collar $\varphi: \partial S_n\times [0,\varepsilon)\to S_n$ which is a representative of the given equivalence class, with $\varepsilon$ small enough such that $V_\alpha\cap S^{(n-1)}\neq \emptyset \Rightarrow \varphi(\hat{f}^{-1}_n ( supp((\lambda_{S^{(n-1)}})_\alpha))\times [0,\varepsilon))\subseteq V_\alpha\cap S_n$. (One can obtain the $\varepsilon$ by going over to a finite subcovering of the $V_\alpha$ at the beginning of the induction process, which is possible since $S$ is cocompact and the $V_\alpha$ are $G$-invariant. Then one can take the minimum of the appearing $\varepsilon$'s.) 

Take a smooth map $f:[0,\varepsilon)\to [0,1]$ such that $f([0,\frac{1}{3}\varepsilon])=0$ and $f([\frac{2}{3}\varepsilon,\varepsilon))=1$, define $f':=fpr_2\varphi^{-1}:\varphi(\partial S_n\times [0,\varepsilon))\to [0,1], \varphi(x,t)\mapsto f(t)$ and continue $f'$ by $1$ to all of $S_n$. Set 
\begin{eqnarray*}
\lambda_{\alpha 1}(x) & := & 
\begin{cases}
f'(x)\cdot (\lambda_{S_n})_\alpha (x) & \text{ if }x\in \oS_n\\
0 & \text{ else}
\end{cases}
\\
\lambda_{\alpha 2}(x) & := &
\begin{cases}
(\lambda_{S^{(n-1)}})_\alpha(x) & \text{ if }x\in S^{(n-1)}\\
\kappa_\alpha(x)
&  \text{ if } x\in \oS_n \text{ and } V_\alpha\cap S^{(n-1)}\neq \emptyset\\
0 &  \text{ else,}
\end{cases}
\end{eqnarray*}
where 
\[
\kappa_\alpha(x)=
\begin{cases}
(1-f')(x)\cdot(\lambda_{S^{(n-1)}})_\alpha(\hat{f_n}pr_1\varphi^{-1}(x)) & \text{ if }x\in \varphi(\hat{f}_n^{-1}(V_\alpha\cap S^{(n-1)})\times (0,\varepsilon))\\
0 & \text{ else, } 
\end{cases}
\]
and define $\lambda_\alpha(x):=\lambda_{\alpha 1}(x) + \lambda_{\alpha 2}(x)$.

The map $\lambda_\alpha$ is a smooth map from the stratifold $S$ to $\R$, since all appearing maps are smooth by construction and $\lambda_\alpha$ fulfills $\lambda_\alpha(\varphi(x,t))=\lambda_\alpha(x)$ for all $x\in \partial S_n$ and $t\in [0,\frac{1}{3}\varepsilon]$. 

The support of $\lambda_\alpha$ is contained in $V_\alpha$. We have $\sum_{\alpha \in A} \lambda_\alpha = 1$, because for $x\in S^{(n-1)}$ we calculate $\sum_{\alpha \in A} \lambda_\alpha (x)= \sum_{\alpha \in A} (\lambda_{S^{(n-1)}})_\alpha (x) = 1$ and for $x\in{\oS}_n$ we calculate
\begin{eqnarray*}
\lefteqn{\sum_{\alpha \in A}\lambda_\alpha(x)}\\
 & = & \sum_{\alpha \in A}(\lambda_{\alpha 1}(x) + \lambda_{\alpha 2}(x))\\
& = & \sum_{\alpha \in A} f'(x)\cdot (\lambda_{S_n})_\alpha (x) + \sum_{\alpha \in A} (1-f')(x)\cdot (\lambda_{S^{(n-1)}})_\alpha(\hat{f_n}pr_1\varphi^{-1}(x)) \\
& = & f'(x)\cdot \sum_{\alpha \in A} (\lambda_{S_n})_\alpha (x) + (1-f')(x) \cdot\sum_{\alpha \in A} (\lambda_{S^{(n-1)}})_\alpha(\hat{f_n}pr_1\varphi^{-1}(x)) \\
& = &  f'(x)\cdot 1 + (1-f')(x)\cdot 1\\
& = & 1.
\end{eqnarray*}
The map $\lambda_\alpha$ is $G$-invariant since all appearing maps are $G$-invariant. So $\{\lambda_\alpha\}$ is the desired smooth $G$-invariant partition of unity on $S$ subordinate to $\{V_\alpha\}_{\alpha\in A}$.
\end{proof}

Note: If we use another definition of collars (with a $\delta$-function), cocompactness is not necessary any more.

\begin{Lem}\label{differentiableGinvariantextension}
Let $M$ be a proper smooth $G$-manifold respectively $S$ a proper $G$-stratifold. Let $A$, $B$ be disjoint closed $G$-invariant subsets of $M$ respectively $S$. Then there is a smooth $G$-invariant map $\rho:M\to \R$ respectively $\rho:S\to \R$ which extends $\rho|_A\equiv 0$, $\rho|_B\equiv 1$. 
\end{Lem}
\begin{proof}
We write out the proof for the case of a manifold $M$ - it is identical for a stratifold $S$, only the existence of a smooth $G$-invariant partition of unity is used. 

The space $G\backslash M$ is a paracompact Hausdorff space, so we can find a continuous map $\rho':G\backslash M\to \R$ extending $\overline{\rho|_{A\cup B}}$. So we also have a continuous $G$-invariant extension $\rho'':M\to \R$ of $\rho|_{A\cup B}$ by setting $\rho''(x):=\rho'(Gx)$ for all $x\in M$. With the smooth $G$-invariant partition of unity on $M$, we construct a smooth $G$-invariant approximation of $\rho'':M\to \R$ also extending $\rho|_{A\cup B}$. (The proof is analogous to~\cite[proof of Th.~II.11.7]{bredon1993}.) For $x\in M$ let $V_x$ be a $G$-invariant neighborhood of $x$ in $M$ and $h_x:V_x\to \R$ such that
\begin{enumerate}
\item $x\in A \Rightarrow V_x\cap B=\emptyset \text{ and }h_x\equiv 0,\; x\in B\Rightarrow  V_x\cap A=\emptyset \text{ and }h_x\equiv 1$. 
\item $x\not\in A\cup B\Rightarrow V_x\cap (A\cup B)=\emptyset$ and $y\in V_x\Rightarrow h_x(y)=\rho''(x)$ (constant in $y$).
\end{enumerate} 
The map $h_x$ is $G$-invariant by construction. (We have $h_{g_1 x}(g_2 y)=\rho''(g_1 x)=\rho''(x)$ for all $g_1, g_2\in G.$)

Let $\{U_\alpha\}$ be a locally finite $G$-invariant refinement of $\{V_x\}$ with index assignment $\alpha\mapsto x(\alpha)$ and let $\{\lambda_\alpha\}$ be a smooth $G$-invariant partition of unity on $M$ with $supp(\lambda_\alpha)\subseteq U_\alpha$. ($\lambda_\alpha=0$ on $A\cup B$ except if $x(\alpha)\in A\cup B$.)

Let $\rho(y):=\sum \lambda_\alpha(y)h_{x(\alpha)}(y)$ for $y\in M$. Then $\rho$ is smooth on $M$ because the $\lambda_\alpha$ and $h_{x(\alpha)}$ are smooth. $\rho$ is $G$-invariant since the $\lambda_\alpha$ and $h_{x(\alpha)}$ are $G$-invariant. If $y\in A$, then $\rho(y)=\sum \lambda_\alpha(y)h_{x(\alpha)}(y)=\sum \lambda_\alpha(y)\cdot 0 = 0$ and if  $y\in B$, then $\rho(y)=\sum \lambda_\alpha(y)h_{x(\alpha)}(y)=\sum \lambda_\alpha(y)\cdot 1 = 1$. So $\rho$ is the desired smooth $G$-invariant extension of $\rho|_A\equiv 0$, $\rho|_B\equiv 1$.  
\end{proof}

Now we are ready to show that our bordism groups form a $G$-homology theory.

\begin{Prop}
Let $\Cl$ be a bordism class of proper $G$-stratifolds. The assignment
\begin{eqnarray*}
\H^G_{\Cl,n}: G-Top & \to & Ab\\
X & \mapsto & \{\; [S,f]\;|\; S \text{ $n$-dim.~cocomp.~proper $G$-stratifold in $\Cl$, }\\
& & \quad \quad \quad \quad f:S\to X \text{ $G$-equivariant}\}\\
(g:X\to Y) & \mapsto & (g_*:[S,f]\mapsto [S,gf])
\end{eqnarray*}
for all $n\in \Nat$ ($\H^G_{\Cl,n}:=0$ for $n< 0$) together with the boundary operator $d$ defined in the proof is a $G$-homology theory. 
\end{Prop}

\begin{proof}
The procedure is analogous to the case without $G$-operation.\\
i) $\H^G_n(\emptyset)=0$ since there is no map of a non-empty stratifold into $\emptyset$.\\
ii) functoriality: clear by definition\\
iii) $G$-homotopy invariance:\\
Let $F:X\times I \to Y$ be a $G$-homotopy of $g_1:X\to Y$ and $g_2:X\to Y$. Let $[V,f]\in \H^G_n(X)$. Then $V\times I$ (with the trivial $G$-operation on $I$) is an $(n+1)$-dimensional cocompact proper $G$-stratifold in $\Cl$ with boundary. We have 
$$ g_{1*} ([V,f]) + g_{2*} ([V,f])=[V,g_1 f]+[V,g_2 f] = [\partial (V\times I),F(f\times id)|_{\partial (V\times I)}] = 0,$$
so $ g_{1*}([V,f])=-g_{2*}([V,f])=g_{2*}([V,f]).$\\
iv) Mayer-Vietoris-sequence:\\
For $X=X_1\cup X_2$, where $X_1$ and $X_2$ are $G$-invariant open subsets of $X$, a boundary operator $d:H^G_n(X_1\cup X_2)\to H^G_{n-1}(X_1\cap X_2)$ is needed (for all $n \in \Nat$). We construct it as follows: If $[S,f]\in H^G_n(X_1\cup X_2)$, then both $f^{-1}(X_1\setminus X_2)$ and $f^{-1}(X_2\setminus X_1)$ are disjoint closed $G$-invariant subspaces of $S$. We define a $G$-invariant map
$$ \tilde{\rho}:f^{-1}(X_1\setminus X_2)\amalg f^{-1}(X_2\setminus X_1)\to \mathbb{R}$$
by setting $$\tilde{\rho}(f^{-1}(X_1\setminus X_2) )=0, \quad\tilde{\rho}(f^{-1}(X_2\setminus X_1) )=1.$$ 
By Lemma~\ref{differentiableGinvariantextension} we can find a differentiable $G$-invariant map $\rho:S\to \mathbb{R}$ extending $\tilde{\rho}$.
%nur f\"ur kompakte Gruppen: The Tietze-Gleason equivariant extension theorem (\cite[I.2.3, together with 0.3.3]{bredon1972}) gives us a differentiable $G$-invariant map $\rho:S\to \mathbb{R}$ extending $\tilde{\rho}$ (by ``averaging out'' the action of $G$, for example for $G$ finite $\rho(x)=1/|G|\sum_{g\in G}\tilde{\rho}'(gx)$, otherwise use the Haar measure on $G$). \\
We choose a regular value $t\in (0,1)$ of $\rho$. We define  $Z:=\rho^{-1}(t)$. This is an $(n-1)$-dimensional stratifold and a $G$-invariant subspace of $S$ (since $\rho$ is $G$-invariant), thus a $G$-stratifold with proper $G$-action (since the action of $G$ on $S$ is proper). It is in $\Cl$ by axiom 3. Since $t\notin \{0,1\}$, we have $Z\not\subseteq f^{-1}(X_1\setminus X_2 \cup X_2 \setminus X_1)$, and so $f(Z)\subseteq X_1\cap X_2$. We set $d([S,f]):=[Z,f|_Z]\in H^G_{n-1}(X_1\cap X_2)$. 

The map $d$ is well defined: Another choice $t'\in (0,1)$ of of the regular value would lead to a singular stratifold $(Z':=\rho^{-1}(t'),f|_{Z'})$ which is bordant to $(Z,f|_{Z})$ by $(\rho^{-1}([t,t']), f|_{\rho^{-1}([t,t'])})$ if $t'>t$ and by  $(\rho^{-1}([t',t]), f|_{\rho^{-1}([t',t])})$ if $t'<t$. Another choice $(S',f')$ of the representative of $[S,f]$ would be bordant to $(S,f)$ by a pair $(T,F)$, and we could define $\rho$ on all of $T$ as above, which would lead to a bordism between $(\rho^{-1}(t)\cap S, f|_{\rho^{-1}(t)\cap S})$ and $(\rho^{-1}(t)\cap S', f|_{\rho^{-1}(t)\cap S'})$ given by $(\rho^{-1}(t), F|_{\rho^{-1}(t)})$ (where $t$ is a regular value of $F$, $f$ and $f'$).

The construction is independent of the choice of $\rho$, which is shown analogously to the non-equivariant case~\cite{kreck}, which in turn is shown analogously to the case of usual bordism of manifolds~\cite{tomdiecktopologie}. We have to show that if there are two $G$-invariant maps $\rho:S\to \R$ and $\rho':S\to \R$ both sending $f^{-1}(X_1\setminus X_2) $ to $0$ and $f^{-1}(X_2\setminus X_1) $ to $1$, then there is a regular value $t$ of $\rho$ and $\rho'$ such that $[\rho^{-1}(t), f|_{\rho^{-1}(t)}]=[\rho'^{-1}(t), f|_{\rho'^{-1}(t)}]$. We take a smooth map $h:[0,1]\to [0,1]$ such that $h([0,\frac{1}{3}])=0$ and $h([\frac{2}{3},1])=1$ and define $\phi:S\times I \to \R, \;(s,x)\mapsto \frac{1}{2}(\rho(s)\cdot (1-h)(x) + \rho'(x)\cdot h(x))$. We know that $S\times I$ is an $(n+1)$-dimensional $G$-stratifold with boundary $S\times\{0\}\amalg S\times\{1\}$, that $\phi$ is a $G$-invariant map, and that $\phi|_{S\times\{0\}\amalg S\times\{1\}}=\rho\amalg \rho'$. We choose a regular value $t$ of $\phi$, $\rho$ and $\rho'$, then $(\phi^{-1}(t),f pr_S|_{\phi^{-1}(t)})$ is a bordism between $(\rho^{-1}(t), f|_{\rho^{-1}(t)})$ and $(\rho'^{-1}(t), f|_{\rho'^{-1}(t)})$. 

The exactness of the Mayer-Vietoris-sequence is also proven analogously to the non-equivariant case~\cite{kreck}, respectively analogously to usual bordism of manifolds~\cite[Satz VIII.13.8]{tomdiecktopologie}. We will not write out the details here, since there are no additional technical difficulties involved.
\end{proof}

We can also place restrictions on the nature of the $G$-operations allowed on the stratifolds. We are especially concerned with ``free proper $G$-actions'' or ``all proper $G$-actions'', but one could imagine many other restrictions such as ``fixed point free'' or, more generally, ``with isotropy in $\F$'', where $\F$ is a set of subgroups of $G$ closed under conjugation or subconjugation. These other restrictions might also lead to interesting homology theories.

\section{Homology theories defined by stratifolds with free proper $G$-Operations} \label{sec2}

If we restrict ourselves to free proper $G$-operations on a certain class $\Cl$ of stratifolds, we obtain an explicit description of the usual $\Cl$-stratifold homology theory of the Borel construction $EG\times_G X$. We show that $\H^{G, free}_{\Cl,*}(X)\cong \H_{\Cl,*}(EG\times_G X),$ where $\H_\Cl$ stands for the homology theory in question. Our standard examples are $Eh$, Euler homology, and $H(-;\Zz)$ or $H(-;\Z)$, singular homology with $\Zz$- or $\Z$-coefficients. The construction is analogous to the one for bordism of smooth manifolds with free $G$-action~\cite[page 50 ff]{conner-floyd}. We define 
\begin{eqnarray*}
\H^{G, free}_{\Cl,n}(X)= & \{ \;[S,f] \;| & S \text{ $n$-dim.~cocomp.~$G$-str.~in $\Cl$ with free proper} \\
& & \text{$G$-action, $f:S\to X$ $G$-equiv.}\}.
\end{eqnarray*}

We want to give a natural transformation 
$$ \H^{G, free}_{\Cl,*}(-)\to \H_{\Cl,*}(EG\times_G -)$$
and check that it is an equivalence of homology theories. 

Let $X\in G-Top$ be a $G$-space and $[S,f]\in \H^{G, free}_{\Cl,n}(X)$. Since $S$ is a free $G$-stratifold, there is up to $G$-homotopy a unique $G$-equivariant map $h:S\to EG$. (Here we consider $EG$ as a left $G$-space, where the left action of $G$ on $EG$ is given by $G\times EG\to EG, (g,e)\mapsto e g^{-1}$. Then $EG\times_G X=G\backslash EG\times X$.) Thus we have a map $(h,f):S\to EG\times X$. Taking $G$-orbits gives a map $\overline{(h,f)}: G\backslash S\to EG\times_G X$. (We have $(h,f)(s)=(h(s),f(s))$ and $(h,f)(gs)=(h(s)g^{-1},gf(s))$, and these are exactly the elements that are identified in the passage from $EG\times X$ to $EG\times_G X$.)

We have to check that $G\backslash S$ (with trivial $G$-action) is again in $\Cl$.

\begin{Lem}\label{orbitraumstratifold}
If $S$ is a $G$-stratifold with free proper $G$-action, then the orbit space $G\backslash S$
again has a stratifold structure given by $(G\backslash S)_i=G\backslash S_i$, $(f_{G\backslash S})_i=\bar{f_i}$. 
\end{Lem} 

\begin{proof}
Since $gS_i=S_i$ for all $g\in G$, every $S_i$ is a free proper $G$-manifold, and thus $G\backslash S_i$ is again a manifold. We have $\partial (G\backslash S_i)=G\backslash \partial S_i$. Collars are also preserved. Since $f_i: \partial S_i\to S^{(i-1)}$ is $G$-equivariant, it induces a map $\bar{f_i}: G\backslash \partial S_i\to G\backslash S^{(i-1)}$ on the orbit spaces. ($\bar{f_i}$ is proper since $f_i$ is proper.)
\end{proof}

\begin{Lem} \label{propertiesfulfilled}
If $S$ is in a class $\Cl$ of stratifolds defined ``locally'' or by emptyness conditions on the strata, $G\backslash S$ is in the same class.
\end{Lem}
\begin{proof} We proceed case by case.
\begin{itemize}
\item $\Cl$ ``locally'' defined: A $G$-invariant neighborhood of $x\in S$ is given by $GV_{U_x}\cong G(B^i\times F)=G B^i\times F$, and this is a neighborhood of $Gx$ in $S$, thus $G\backslash GV_{U_x}\cong G\backslash (G B^i\times F)\cong B^i\times F$ is a neighborhood of $Gx$ in $G\backslash S$. $F$ remains unchanged. 
\item $\Cl$ defined by emptyness conditions on strata: $(G\backslash S)_i=G\backslash S_i=\emptyset \Leftrightarrow S_i=\emptyset$, so properties remain unchanged.
\item $\Cl$ with orientability condition: If $S_n$ is orientable and $G$ operates orientation preserving, then $G\backslash S_n$ is orientable.
\end{itemize}
\end{proof}

We know that $G\backslash S$ is compact since $S$ is cocompact. Thus we have an assignment
\[
\H^{G, free}_{\Cl,n}(X) \to \H_{\Cl,n}(EG\times_G X), \quad
{}[S,f] \mapsto [G\backslash S, \overline{(h,f)}].
\]

\begin{Prop}
The assignment 
\begin{eqnarray*}
\H^{G, free}_{\Cl,n}(X) & \to & \H_{\Cl,n}(EG\times_G X)\\
{}[S,f] & \mapsto & [G\backslash S, \overline{(h,f)}].
\end{eqnarray*}
induces a natural equivalence of $G$-homology theories on the category $G-Top$.
\end{Prop}
\begin{proof} We proceed in several steps.
\begin{itemize}
\item The assignment is well defined: if we choose another representative $(S',f')\in [S,f]$, then there is a free proper $G$-stratifold $T$, $\partial T=S\amalg S'$, $F:T\to X$, $F|_{\partial T}=f\amalg f'$. So $(G\backslash T, \overline{(h,F)})$ gives a bordism between $(G\backslash S, \overline{(h,f)})$ and $(G\backslash S', \overline{(h,f')})$.
\item This assignment is evidently a group homomorphism.
\item It is functorial in $X$: Let $g:X\to Y$ be a $G$-equivariant map. Then $$
\xymatrix{
\H^{G, free}_{\Cl,n}(X) \ar[r] \ar[d]^{g_*} & \H_{\Cl,n}(EG\times_G X) \ar[d]^{(id\times g)_*}\\
\H^{G, free}_{\Cl,n}(Y) \ar[r] & \H_{\Cl,n}(EG\times_G Y)
}$$
commutes since $(id\times g)_* \overline{(h,f)} = \overline{(h,gf)}:S\to EG\times_G Y$.
\item It is compatible with the boundary operator since we chose $\rho$ $G$-invariant. (We see that a choice of $\rho$ for $\H^{G, free}_{\Cl,n}(X)$ induces a choice of $\rho$ for $\H_{\Cl,n}(EG\times_G X)$ and vice versa, since we otherwise made the same construction as for usual stratifold homology theories. Note that for $G$-invariant $X_1, X_2$ we have $EG\times_G(X_1\cup X_2) = EG\times_G X_1 \cup EG\times_G X_2$.)
\item It is an isomorphism for each $X\in G-Top$, since we can construct an inverse $\H_{\Cl,n}(EG\times_G X) \to  \H^{G, free}_{\Cl,n}(X)$: Let $[W,j]\in \H_{\Cl,n}(EG\times_G X)$. We have the following commutative diagram:$$
\xymatrix{
S:=(pr_1 j)^* EG \ar[rr]^h \ar[d]^\pi & & EG \ar[d]^\nu \\
W \ar[r]^j & EG\times_G X \ar[r]^{pr_1} & BG 
}$$
$S$ is a $G$-principal bundle over $W$. $S$ is a stratifold: Consider the diagram$$
\xymatrix{
S_i:=(pr_1 j \; inc_i)^* EG \ar[dd]^{\pi_i} \ar@{..>}[dr] \ar[drrr]& & &\\
& S \ar[rr]^h \ar[dd]^\pi & & EG \ar[dd]^\nu \\
W_i \ar[dr]^{inc_i} &&& \\
& W \ar[rr]^{pr_1 j} & & BG 
}$$
This defines strata $S_i$ and unique $G$-equivariant inclusion maps into $S$ via the pullback property. $S_i$ is an $i$-dimensional smooth manifold since $\pi_i:S_i\to W_i$ is a local homeomorphism. For the same reason collars are preserved. $S_i$ is a covering space of $W_i$ with fiber $G$ (discrete). The free action of $G$ on the fiber induces a free action of $G$ on $S_i$. $S$ is a covering space of $W$ with fibre $G$, so $G\backslash S \cong W$ is compact and we have a free action of $G$ on $S$. Since $S\to W$ is a local homeomorphism, $S\in \Cl \Leftrightarrow W\in \Cl $ for stratifold classes $\Cl$ defined locally or by emptyness conditions. 

It is known that $$
\xymatrix{
EG\times X \ar[d] \ar[r]   & EG \ar[d] \\
EG\times_G X  \ar[r] & BG 
}$$
is a pullback diagram (\cite[II.2.5]{bredon1972}). Thus there is a unique $G$-equivariant map $\tilde{f}$ making the following diagram commutative:$$
\xymatrix{
S \ar[dd]^{\pi} \ar@{..>}[dr]_{\tilde{f}} \ar[drrr]^h & & &\\
& EG\times X \ar[dd] \ar[rr]   & & EG \ar[dd] \\
W \ar[dr]^j &&& \\
& EG\times_G X  \ar[rr]^{pr_1} & & BG 
}$$
If we take this map $\tilde{f}$ and compose it with the projection $pr_X$ to $X$, we obtain a map $f:=pr_X \tilde{f}:S\to X$. Thus we have constructed a pair $(S,f)$ such that $[S,f]\in\H^{G,free}_{\Cl,n}(X)$ and which maps to $[W,j]$ under the above correspondence: $G\backslash S=W$, $\overline{(h,f)}=\overline{\tilde{f}} =j$. 
\end{itemize}
\end{proof}

\Rem
If we want to extend this method to other restrictions on the $G$-actions by replacing $EG$ by an appropriate $E(G,\F)$, we have to think about whether $G\backslash S$ is still a stratifold and still in $\Cl$. We might make it into a stratifold by giving it a more complicated stratifold structure, but then it is in general not in $\Cl$ any more.

\section{The induction structure} \label{sec3}

From now on, we restrict ourselves to $G$-CW-complexes $X$: In the sequel, we mean functors from $G$-CW-complexes to $Ab$ when we speak about $G$-homology theories.

Let $\Cl$ be a bordism class such that 
\begin{itemize}
\item If $S\in\Cl$ and $K$ acts freely and properly on $S$ ($K$-cocompact), then $K\backslash S$ is in $\Cl$.
\item If $S\in\Cl$, then any connected component of $S$ is in $\Cl$ again.
\item If $S\in\Cl$, then a discrete union of copies of $S$ is in $\Cl$.
\end{itemize} 
The fact that $K\backslash S$ is again a stratifold was proven in Lemma~\ref{orbitraumstratifold}. The above properties are fulfilled for the usual restrictions: The first point was proven in Lemma~\ref{propertiesfulfilled}, the other points are shown completely analogous.

We now show that we have an induction structure on our $G$-homology theories (as defined for example in~\cite[page 198f]{lueck}). This induction structure links the various homology theories for different groups $G$. 

%wollen wir G-CW-Komplexe oder propere G-CW-Komplexe?????

%$ker (\alpha)$ acts freely$\to$ properly, da wir G-CW-Komplexe haben 
\begin{Prop}
We have an induction strucure on the $G$-homology theories 
$$\H^G_{\Cl,*}: G\text{-CW-complexes }\to Ab$$ 
defined as above by 
$$\H^G_{\Cl,n}(X)= \{ [S,f]\,|\,S \text{ $n$-dim.~cocomp.~pr.~$G$-str.~in $\Cl$, }f:S\to X \text{ $G$-equiv.}\}.$$
Given a group homomorphism $\alpha: H\to G$ and an $H$-CW-complex $X$ such that $ker (\alpha)$ acts freely on $X$, we define maps
\begin{eqnarray*}
ind_\alpha : \H^H_{\Cl,n}(X) & \to & \H^G_{\Cl,n}(ind_\alpha (X))\\
{}[S,f] & \mapsto & [ind_\alpha (S), ind_\alpha(f)].
\end{eqnarray*}
Here $ind_\alpha (X)=G\times_\alpha X$, $ind_\alpha (S)=G\times_\alpha S$ and $ind_\alpha(f)=id_G\times_\alpha f$.

These maps $ind_\alpha$ are isomorphisms for all $n\in \Z$ and satisfy 

\noindent{\bf a) compatibility with boundary homomorphisms}
$$ d_n^G\; ind_\alpha=ind_\alpha\; d_n^H.$$

\noindent{\bf b) functoriality}

Let $\beta:G\to K$ be another group homomorphism such that $ker(\beta\alpha)$ acts freely on $X$. Then we have 
$$ind_{\beta\alpha}=\H_n^K(f_1)ind_\beta ind_\alpha:\H_n^H(X)\to\H_n^K(ind_{\beta\alpha}(X)),$$
where $f_1:ind_\beta ind_\alpha(X)\xrightarrow{\cong} ind_{\beta\alpha}(X),\quad [k,g,x]\mapsto [k\beta(g),x]$ is the natural $K$-homeomorphism.

\noindent{\bf c) compatibility with conjugation}

For $g\in G$ and a $G$-CW-complex $X$ the homomorphism
$$ind_{c(g):G\to G}:\H_n^G(X)\to\H_n^G(ind_{c(g):G\to G}(X))$$
agrees with $\H_n^G(f_2)$ for the $G$-homeomorphism $f_2:X\to ind_{c(g):G\to G}(X)$ which sends $x$ to $[1,g^{-1} x]$ in $G\times_{c(g)} X$.
\end{Prop}

\begin{proof}
It is clear that we only have to treat the case $n\geq 0$ since for $n<0$ all appearing homology groups are $0$. We first need to show that $ind_\alpha$ is well defined:

Let $S$ be an $n$-dimensional proper cocompact $H$-stratifold in $\Cl$. Then $ind_\alpha(S)=G\times_\alpha S=(\cup_{\bar{g}\in G / \alpha(H)}\bar{g} \alpha(H))\times_\alpha S=\cup_{\bar{g}\in G / \alpha(H)}\bar{g} (\alpha(H)\times_\alpha S)$, and $\alpha(H)\times_\alpha S\cong ker(\alpha)\backslash S$. Note that $ker (\alpha)$ acts freely on $S$ since $f$ is $H$-equivariant and $ker (\alpha)$ acts freely on $X$. It acts properly on $S$ since $H$ acts properly on $S$. So $ker (\alpha)\backslash S$ is in $\Cl$ again. It is an $n$-dimensional stratifold since $ker(\alpha)$ is a discrete group. A discrete union of copies of this is in $\Cl$ again. Thus $ind_\alpha(S)$ is in $\Cl$ again. The $n$-dimensional stratifold $ind_\alpha(S)$ is endowed with the obvious $G$-action by left multiplication in the first variable. It is a cocompact $G$-stratifold, since $G\backslash ind_\alpha(S)=G\backslash G\times_\alpha S \cong H \backslash S$ is compact. 

It is a proper $G$-stratifold: We know that a discrete group $G$ acts properly on a Hausdorff space $X$ if and only if for each pair of points $(a,b)$ in $X$ there exist neighborhoods $V_a$ of $a$ and $V_b$ of $b$ such that the set $\{g\in G\;|\;gV_a\cap V_b\neq \emptyset\}$ is finite~\cite[Cor.~I.3.22]{tomdieck}. Thus we need to show that for all $[g,x]$ and $[g',x']$ in $G\times_\alpha S=ind_\alpha(S)$ there are open neighborhoods $W_x$ of $[g,x]$ and $W_{x'}$ of $[g',x']$ such that the set $\{\tilde{g}\;|\; \tilde{g} W_x\cap W_{x'}\neq \emptyset\}$ is finite. We use that $S$ is a proper Hausdorff $H$-space with $H$ discrete, thus for all $x,x'\in S$ there are $V_x,V_{x'} \text{ open such that } \{h\in H\;|\;hV_x\cap V_{x'}\neq \emptyset\}=:\tilde{H}$ is finite. We set $W_x:=\overline{g\times V_x}$ and $W_{x'}:=\overline{g'\times V_{x'}}$ and claim that $\{ \tilde{g}\;|\;\tilde{g}W_x\cap W_{x'}\neq \emptyset \}=g'\alpha(\tilde{H})g^{-1}$, thus finite. This would show that $ind_\alpha (S)$ is a proper $G$-space. 

Proof of the claim:

We have $\tilde{g}W_x\cap W_{x'}\neq \emptyset$ if and only if there are $[\lambda_1,a_1]\in W_x$ and $[\lambda_2,a_s]\in W_{x'}$ such that $[\tilde{g}\lambda_1,a_1]=[\lambda_2,a_2].$
Since ${}[\lambda_1,a_1]\in W_x$, there is $h'_1\in H$ such that $\lambda_1\alpha(h'^{-1}_1)=g \text{ and } h_1'a_1\in V_x.$
Since ${}[\lambda_2,a_2]\in W_{x'}$, there is $ h'_2\in H$ such that $\lambda_2\alpha(h'^{-1}_2)=g'\text{ and } h_2'a_2\in V_{x'}.$ Because of 
${}[\tilde{g}\lambda_1,a_1]=[\lambda_2,a_2]$ there are $h_1,h_2\in H \text{ such that }\tilde{g}\lambda_1\alpha(h_1)=\lambda_2\alpha(h_2)$ and $h_1^{-1}a_1=h_2^{-1}a_2.$

Thus we obtain the equation $\tilde{g}g\alpha(h_1')\alpha(h_1)=g'\alpha(h_2')\alpha(h_2)$ which implies $\tilde{g}=g'\alpha(h_2'h_2h_1^{-1}h'^{-1}_1)g^{-1},$ and $h_2h_1^{-1}a_1=a_2$ which in turn implies that $h_2'h_2h_1^{-1}h'^{-1}_1(h_1'a_1)=(h_2'a_2)$. So $h_2'h_2h_1^{-1}h'^{-1}_1=:h\in\tilde{H}$, and we know that $\tilde{g}=g'\alpha(h)g^{-1}\in g'\alpha(\tilde{H})g^{-1}$.

Conversely, if $\tilde{g}\in g'\alpha(\tilde{H})g^{-1}$, then $\tilde{g}=g'\alpha(h)g^{-1}$ for an $h\in\tilde{H}$, and there are $a_1\in V_{x}$ and $a_2\in V_{x'}$ such that $ha_1=a_2$. We compute $\tilde{g}[g,a_1]=[g'\alpha(h)g^{-1}g,a_1]=[g',ha_1]=[g',a_2]$, thus $\tilde{g}W_x\cap W_{x'}\neq \emptyset.$

 So $ind_\alpha(S)$ is an $n$-dimensional proper cocompact $G$-stratifold in $\Cl$ and $ind_\alpha(f) : ind_\alpha(S)\to ind_\alpha(X)$ is a $G$-equivariant continuous map by definition, thus the map $ind_\alpha$ is well defined.

{\bf Step 1:} The induction map $ind_\alpha : \H^H_{\Cl,n}(X) \to \H^G_{\Cl,n}(G\times_\alpha X)$ is an isomorphism.

We want to define an inverse map $\psi:\H^G_{\Cl,n}(G\times_\alpha X) \to \H^H_{\Cl,n}(X)$. Let $[S',f']\in \H^G_{\Cl,n} (G\times_\alpha X)$, so $S'$ is an $n$-dimensional cocompact proper $G$-stratifold and $f':S'\to G\times_\alpha X$ is a $G$-equivariant map. We have $\alpha(H)\leq G$, and so $\alpha(H)\times_\alpha X \subseteq G\times_\alpha X$, which is an open and closed subspace since $G$ is discrete. 

The operation of $G$ on $G\times_\alpha X$ induces an operation of $H\stackrel{\alpha}{\hookrightarrow} G$ on $G\times_\alpha X$ by $h.[g,x]:=[\alpha(h)g,x]$. The subspace $\alpha(H)\times_\alpha X$ is invariant under this operation.

Define $W:=f'^{-1}(\alpha(H)\times_\alpha X)\subseteq S'$. The map $f'$ is continuous, so $W$ is open and closed in $S'$. Thus $W$ is an $n$-dimensional stratifold in $\Cl$ again. 

$H$ operates on $W$ by $H\times W \to W, (h,w)\mapsto \alpha(h)w$ since $f'$ is $G$-equivariant. This is a restriction of a proper $G$-action, thus proper. Note that $G\times_{\alpha} W \cong S'$ since $f'$ is $G$-equivariant. Thus $W$ is a cocompact $H$-space, since $H\backslash W\cong G\backslash G\times_{\alpha}W\cong G\backslash S'$, and $G\backslash S'$ is compact since $S'$ is a cocompact $G$-space. The restriction of $f'$ to $W$ is an $H$-equivariant map $f'|_W: W\to \alpha(H)\times_\alpha X\cong ker(\alpha)\backslash X$. 

Since $ker(\alpha)$ operates freely on $X$, we can lift this: $X\to ker(\alpha)\backslash X$ is a $ker(\alpha)$-principal bundle, and we define $\tilde{W}:=(f'|_W)^*(X)$ as the pullback of this bundle along $f'|_W$:
%operiert automatisch proper, da X CW-Komplex
\[
\xymatrix{
\tilde{W}:=(f'|_W)^*(X) \ar[r]^{\tilde{f}} \ar[d] & X \ar[d]\\
W \ar[r]^{f'|_W} & ker(\alpha)\backslash X
}\]
Concretely, $\tilde{W}=\{(w,x)\in W\times X\;|\;f'(w)=[1,x]\in \alpha(H)\times_\alpha X\}$. The operation of $H$ on $\tilde{W}$ is induced by the action on $W$ and $X$: We have $h.(w,x):=(\alpha(h)w,hx).$ Thus $\tilde{f}$ is $H$-equivariant by definition. $\tilde{W}$ is a cocompact $H$-space: $H\backslash \tilde{W}=\{(Hw,Hf'(w))\subseteq H\backslash W\times H\backslash X\}\cong H\backslash W$, and this is compact.

$H$ operates properly on $\tilde{W}$: $\tilde{W}$ is a Hausdorff space, so $H$ operates properly if and only if for every $(w,x)$ and $(w',x')$ in $\tilde{W}$ there are open neighborhoods $V_{(w,x)}$ and $V_{(w',x')}$ such that the set $\{h\in H\;|\;h V_{(w,x)}\cap V_{(w',x')}\neq \emptyset\}=:\tilde{H}_{\tilde{W}}$ is finite. We know that $H$ operates properly on $W$, so we know that for $w$ and $w'$ in $W$ there are open neighborhoods $V_w$ and $V_{w'}$ such that the set $\{h\in H\;|\; h V_w\cap V_{w'}\neq \emptyset\}=:\tilde{H}_W$ is finite. We set $V_{(w,x)}:=(V_w\times X)\cap \tilde{W}$ and $V_{(w',x')}=(V_{w'}\times X)\cap \tilde{W}$. Then  $\tilde{H}_{\tilde{W}}$ is contained in $\tilde{H}_W$, thus finite.
 
The space $\tilde{W}$ is a $ker(\alpha)$-principal bundle over $W$, thus a stratifold again: The strata are given by the pullbacks of the strata of $W$, which map to $\tilde{W}$ by the universal property of $\tilde{W}$: 
$$
\xymatrix{
\tilde{W}_i:=(f'|_W \; inc_i)^* (X) \ar[dd] \ar@{-->}[dr] \ar[drr]\\
& \tilde{W} \ar[r]_{\tilde{f}} \ar[dd] & X \ar[dd]\\
W_i \ar[dr]^{inc_i} & & \\
& W \ar[r]^{f'|_W} & ker(\alpha)\backslash X
}$$

The $W_i$ are c-manifolds, thus the $\tilde{W_i}$ are c-manifolds, too, since the projection map is a local homeomorphism. The stratifold $\tilde{W}$ is locally isomorphic to $W$. Thus the collars and all usual (local) restrictions are preserved, and $\tilde{W}\in \Cl$ again. We define the map 
\begin{eqnarray*}
\psi: \H^G_{\Cl,n}(G\times_\alpha X)& \to & \H^H_{\Cl,n}(X) \\
{}[S',f'] & \mapsto & [\tilde{W}, \tilde{f}].
\end{eqnarray*}
We now show that $\psi$ is inverse to $ind_\alpha$. 

{\bf 1.} The composition 
\begin{eqnarray*}
 ind_\alpha \psi :  \H^G_{\Cl,n}(G\times_\alpha X) & \to & \H^G_{\Cl,n}(G\times_\alpha X) \\
{}[S',f'] & \mapsto [\tilde{W}, \tilde{f}] \mapsto & [G\times_\alpha \tilde{W}, id_G\times_\alpha \tilde{f}]
\end{eqnarray*}  
 is equal to the identity on $\H^G_{\Cl,n}(G\times_\alpha X)$: We show that 
$$
\xymatrix{
G\times_\alpha \tilde{W} \ar[rr]_{id_G\cdot pr_W}^{\simeq} \ar[dr]_{id_G\times \tilde{f}} & & S' \ar[dl]^{f'}\\
& G\times_\alpha X &
}$$
commutes, where all appearing maps are $G$-equivariant  and the horizontal map $id_G\cdot pr_W$ is an isomorphism of $G$-stratifolds.

The map $f'$ is $G$-equivariant by definition. The map $id_G\times \tilde{f}$ is $G$-equivariant:
$(id_G\times\tilde{f})(\tilde{g}[g,(w,x)])=[\tilde{g}g,\tilde{f}((w,x))]=\tilde{g}[g,x]=\tilde{g}(id_G\times\tilde{f}([g,(w,x)])).$ The map $id_G\cdot pr_W$ is $G$-equivariant: $(id_G\cdot pr_W)(\tilde{g}[g,(w,x)])=\tilde{g}g.w=\tilde{g}(id_G\cdot pr_W)([g,(w,x)]).$

We define an inverse $\varphi: S'\to G\times_\alpha \tilde{W}$ to $id_G\cdot pr_W$: If $s'\in S'$, then $f'(s')= [g,x]\in G\times_\alpha X$, and we set $\varphi(s'):=[g,(g^{-1}s', x)] \in G\times_{\alpha} \tilde{W}$. (Remember that $\tilde{W}=\{(w,x)\;|\;f'(w)=[1,x]\in \alpha(H)\times_\alpha X\}\subseteq W\times X$ and $W=f'^{-1}(\alpha(H)\times_\alpha X)\subseteq S'$.) This map $\varphi$ is well defined: We know that $g^{-1}s'\in W$ since $f'(g^{-1}s')=g^{-1}f'(s')=g^{-1}[g,x]=[1,x]\in\alpha(H)\times_\alpha X.$ 

If we had chosen a different representative $[g\alpha(h^{-1}), hx]=f'(s')\in G\times_\alpha X$, this would lead to the same element: 
$[g\alpha (h^{-1}), ((g\alpha(h^{-1}))^{-1}s',hx)] = [g\alpha (h^{-1}),(\alpha(h) g^{-1}s',hx)]= [g\alpha (h^{-1}),h.(g^{-1}s',x)]=\varphi(s').$

The map $\varphi$ is $G$-equivariant: What is $\varphi(\tilde{g}s')$? We know that $\tilde{f}(\tilde{g} s')=\tilde{g} \tilde{f}(s')=\tilde{g}[g,x]=[\tilde{g}g,x].$ Thus we have $\varphi(\tilde{g}s')=[\tilde{g}g, ((\tilde{g}g)^{-1}(\tilde{g}s'),x)]=[\tilde{g}g, (g^{-1}\tilde{g}^{-1}\tilde{g}s',x)]=\tilde{g}[g,(g^{-1}s',x)]=\tilde{g}\varphi(s').$

The map $\varphi$ is indeed an inverse to $(id_G\cdot pr_W)$: We calculate
$$ (id_G\cdot pr_W)\varphi(s')=(id_G\cdot pr_W)[g,(g^{-1}s', x)]=gg^{-1}s'=s'. $$
Conversely, we calculate
$\varphi(id_G\cdot pr_W)[g,(w,x)]=\varphi(gw).$ We know $f'(gw)=gf'(w)=g [1,x] = [g,x],$
so by definition $\varphi(gw) = [g, (g^{-1}(gw),x)]=[g,(w,x)].$
Both $(id_G\cdot pr_W)$ and $\varphi$ preserve the stratifold structure, thus this is an isomorphism of $G$-stratifolds.

The diagram commutes: If $f'(s')=[g,x]$, then $(id_G\times \tilde{f})\varphi (s')=(id_G\times \tilde{f})[g,(g^{-1}s',x)]=[g,x].$

{\bf 2. }The composition 
\begin{eqnarray*}
\psi \; ind_\alpha: \H^H_{\Cl,n}(X) & \to & \H^H_{\Cl,n}(X), \\
{}[S,f] & \mapsto [\underbrace{ind_\alpha(S)}_{S'}, \underbrace{ind_\alpha(f)}_{f'}] \mapsto & [\tilde{W}, \tilde{f}]
\end{eqnarray*}  
 is equal to the identity on $\H^H_{\Cl,n}(X)$ since  $$
\xymatrix{
S \ar[rr]^{ s\mapsto ([1,s], f(s))}_{\simeq} \ar[dr]_{f}& & \tilde{W} \ar[dl]^{\tilde{f}}\\
& X & 
}$$
commutes, all appearing maps are $H$-equivariant, and the horizontal map is an isomorphism of $H$-stratifolds. (Remember that $\tilde{W}=\{(w,x)\in W\times X | f'(w) =[1,x] \in G\times_\alpha X\}$ and $W=f'^{-1}(\alpha(H)\times_\alpha X)=\alpha(H)\times_\alpha f^{-1}(X)=\alpha(H)\times_\alpha S$.)

The horizontal map is injective: $([1,s],f(s))=([1,s'], f(s'))$ implies $[1,s]=[1,s']$ and $f(s)=f(s')$. This implies $s'=hs$ with $h\in ker(\alpha)$, so $f(s)=hf(s)$. But $ker(\alpha)$ operates freely on $X$, so $h=1$, thus $s'=s$.

It is surjective: If $(w,x)\in \tilde{W}$, then $w=[1,s]\in W$, and for $f(s)\in X$ we have $[1,f(s)]=[1,x]\in G\times_\alpha X$. This means there is an $h\in ker(\alpha)$ such that $hf(s)=x$, or equivalently $f(hs)=x$. We see that $hs$ is mapped to $(w,x)$: $hs \mapsto ([1,hs],f(hs))=([\alpha(h),s],x)=([1,s],x)=(w,x).$

It is $H$-equivariant: For $\tilde{h}\in H$ we obtain the map $\tilde{h}s\mapsto ([1,\tilde{h}s],f(\tilde{h}s))=([\alpha(\tilde{h}),s],\tilde{h} f(s))=\tilde{h}.([1,s],f(s)).$ The horizontal map $S\to \tilde{W}$ is open: The projection $\tilde{W}\to W$ is a local homeomorphism, so $S\to \tilde{W}$ is open if $S\to\tilde{W}\to W$ is open. But this is the map $s\mapsto([1,s],f(s))\mapsto[1,s]$, which is open since $W$ has the quotient topology. The map obviously respects the stratifold structure, so it is an isomorphism of $G$-stratifolds. The diagram clearly commutes: $\tilde{f}([1,s],f(s))=f(s)$.

{\bf Step 2: }The map $ind_\alpha$ satisfies properties a), b) and c):

a) 
We check the compatibility with the boundary operator in the Mayer-Vietoris-sequence. Let $X=X_1\cup X_2$ be an $H$-space, with $X_1$ and $X_2$ $H$-invariant. Note that $ind_\alpha (X_1\cup X_2)=ind_\alpha (X_1) \cup ind_\alpha (X_2)$ and $ind_\alpha (X_1\cap X_2)=ind_\alpha (X_1) \cap ind_\alpha (X_2)$ because $X_1$ and $X_2$ are $H$-invariant. The diagram $$
\xymatrix{
H^H_n(X_1\cup X_2) \ar[rr]^{ind_\alpha} \ar[dd]^{d_H} & & H^G_n(ind_\alpha (X_1) \cup ind_\alpha (X_2)) \ar[dd]^{d_G}\\
&&\\
H^H_{n-1}(X_1\cap X_2) \ar[rr]^{ind_\alpha} & & H^G_{n-1}(ind_\alpha (X_1) \cap ind_\alpha (X_2))
}$$
commutes since we can choose the $G$-invariant map $\rho_G$ on the right (used for the definition of $d_G$) to be $ind_\alpha(\rho_H)$, with $\rho_H$ the $H$-invariant map on the left side. (The construction is independent of the choice of $\rho$.)

b) The diagram \[
\xymatrix{
\H^H_n(X) \ar[r]^-{ind_\alpha}
         \ar`u[rrr]`[rrr]^{ind_{\beta\alpha}}
& \H^G_n(G\times _\alpha X) \ar[r]^-{ind_\beta}
& \H^K_n(K\times_\beta G \times_\alpha X) \ar[r]^-{\H^K_n(f_1)}
& \H^K_n(K\times_{\beta\alpha} X),
}\]
%  \[
%  \xymatrix{
%  \Sigma^m E_n \ar[d] \ar`l[dd]`[dd]_-{\sigma}[dd] \ar[r]^-{\Sigma^m f_n} &
%  \Sigma^m F_n \ar[d] \ar`r[dd]`[dd]^-{\sigma}[dd] \\
%  E_{(m,n)} \ar[d] \ar[r]^-{f_{(m,n)}} &
%  F_{(m,n)} \ar[d] \\
%  E_{m+n} \ar[r]^-{f_{m+n}} & F_{m+n}
%  }
% \]
where $f_1=id_K\cdot \beta \times_\alpha id_X: K\times_\beta G \times_\alpha X\to K\times_{\beta\alpha} X$ is the natural $K$-homeomorphism,   
commutes if $[K\times_\beta G\times_\alpha S, f_1(id_K\times_\beta id_G\times_\alpha f)]=[K\times_{\beta\alpha} S, id_K\times_{\beta\alpha} f].$ This is true since the following diagram 
$$
\xymatrix{
K\times_\beta G\times_\alpha S \ar[dr]_{f_1(id_K\times_\beta id_G\times_\alpha f)} \ar[rr]^{\simeq}_{id_K\cdot \beta \times_\alpha id_S} & & K\times_{\beta\alpha} S \ar[dl]^{id_K\times_{\beta\alpha} f}\\
& K\times_{\beta\alpha} X & 
}$$
commutes, all appearing maps are $K$-equivariant,  and the horizontal map is an isomorphism of $K$-stratifolds. The $K$-equivariance of the maps is immediate. The horizontal map is the natural $K$-homeomorphism $ind_\beta ind_\alpha(S)\xrightarrow{\sim} ind_{\beta\alpha}(S), (k,g,s)\mapsto(k\beta(g),s).$ It obviously respects the stratifold structure. The commutativity of the diagram is clear: We calculate $(id_K\times_{\beta\alpha}f)(id_K\cdot\beta\times_\alpha id_S)(k,g,s)=(k\beta(g),f(s))=f_1(id_K\times_\beta id_G\times_\alpha f)(k,g,s).$\\

c) The diagram $$ 
\xymatrix{
\H^G_n(X) \ar@/_/[d]_{\H_n(f_2)} \ar@/^/[d]^{ind_c(g)}\\
\H^G_n(G\times_{c(g)} X),
}$$
where $f_2:X\to G\times_{c(g)}X, x\mapsto (1,g^{-1}x)$ is a $G$-homeomorphism, commutes since the following diagram commutes, all appearing maps are $G$-equivariant, and the horizontal map is an isomorphism of $G$-stratifolds:$$
\xymatrix{
S \ar[dr]_{f_2 f} \ar[rr]^{\simeq}_{s\mapsto[1,g^{-1}s]} & & G\times_{c(g)} S \ar[dl]^{id_G\times_{c(g)}f}\\
& G\times_{c(g)} X &
}$$
The horizontal map is the $G$-homeomorphism $S\to G\times_{c(g)} S, s\mapsto [1,g^{-1}s]$ which obviously respects the stratifolds structure. The maps are $G$-equivariant since $f_2$ is $G$-equivariant. The commutativity of the diagram is clear: We calculate $(id_G\times_{c(g)} f)(1,g^{-1}s)=[1,f(g^{-1}s)]=[1,g^{-1}f(s)]=f_2 f(s).$
\end{proof}

We have shown that our construction yields $G$-homology theories with an induction structure, thus equivariant homology theories. In the next sections we study some of the homology theories thus obtained in detail.

\section{Homology theories defined by stratifolds with arbitrary proper $G$-Operations} \label{sec4}

% hier brauchen wir proper nicht, da wir nicht zu den Orbitraeumen uebergehen, denke ich
% oben brauchen wir proper schon, denn da sind wir zu S/G uebergegangen
% Achtung, es gibt viele Sachen, die nicht gehen, wenn die Geschichte nicht proper ist: lernen!!!

We use the $G$-comparison theorem and the cone construction, refining the procedure employed in the analysis of usual stratifold homology theories. In order to use the $G$-comparison theorem, we need to compute the coefficients for all homogeneous $G$-spaces, thus for all $G/H$ with $H\leq G$. The induction map applied to the inclusion $\alpha: H\hookrightarrow G$ gives an isomorphism
$$ ind_\alpha: \H^H_{\Cl,n}(*)\to \H^G_{\Cl,n}(\underbrace{ind_\alpha *}_{G/H}).$$
This allows us to reduce our analysis of $\H^G_{\Cl,n}$ to the inspection of $\H^H_{\Cl,n}(*)$ for all $H\leq G$. 

The cone construction only works for actions of finite groups $H$, since otherwise the trivial action on the cone point is not proper any more. So we only obtain complete answers in the case of finite groups $G$ or proper $G$-CW-complexes $X$.

Let $H$ be a finite group. By forgetting the $H$-action, we get a map $\H^H_{\Cl,n}(*)\to \H_{\Cl,n}(*)$. This is well-defined since the finiteness of $H$ implies that an $H$-cocompact stratifold is compact. The map is split surjective since we can endow each stratifold $S$ with the trivial $H$-action (which is proper since $H$ is finite). Thus we get a refinement of usual homology theories defined via stratifolds. 

Let us now look in detail at the equivariant homology theories defined by different bordism classes $\Cl$.

\subsection{$\Cl=All$, ``all stratifolds''}

\begin{Prop} Let $H$ be a finite group. We obtain $\H^H_{All,n}(*)=0$ for all $n\geq 0$. 
\end{Prop}
\begin{proof}
Let $[S]\in \H^H_{All,n}(*).$ We can ``cone off'' $S$ and extend the given $H$-action to the cone by setting $CS=S\times [0,1]/S\times \{1\}$, $h(x,t):=(hx,t).$ Then $[S]=[\partial CS]=0.$ 
\end{proof}

\begin{Cor}
Let $G$ be a finite group. Then $\H^G_{All,n}(G/H)=0$ for all $H\leq G$ and for all $n$, thus $\H^G_{All,*}=0$ on all $G$-CW-complexes $X$. 
\end{Cor}
\begin{proof}
The induction structure gives us $\H^G_{All,n}(G/H)=\H^H_{All,n}(*)=0$, thus the coefficients are all trivial. Since there is a map to the zero homology theory, we can use the comparison theorem to obtain the fact that the homology theory we have constructed is in fact the zero homology theory. 
\end{proof}

\begin{Cor} \label{homgroupszero}
Let $G$ be a discrete group, and let $H$ be a finite subgroup of $G$. Then $\H^G_{All,n}(G/H)=0$ for all $n\geq 0$. 
\end{Cor}
\begin{proof}
$\H^G_{All,n}(G/H)=\H^H_{All,n}(*)=0.$
\end{proof}

\begin{Cor}
Let $G$ be a discrete group, and let $X$ be a proper $G$-CW-complex. Then $\H^G_{All,*}(X)=0$.
\end{Cor}
\begin{proof}
A proper $G$-CW-complex $X$ is built up of cells of the form $G/H\times D^i$ with $H$ a finite subgroup of $G$, $i\geq 0$. So one can use the Mayer-Vietoris-sequence to reduce the calculation of $\H^G_{All,*}(X)$ to the calculation of $\H^G_{All,n}(G/H)$ for all finite subgroups $H$ of $G$, for $n\geq 0$. But these homology groups are all zero by Corollary~\ref{homgroupszero}. 
\end{proof}

\subsection{$\Cl=$ ``Euler stratifolds''}

Recall the following definition:
\begin{Def}
An $n$-dimensional stratifold $S$ is called \emph{regular} if for all $0\leq i \leq n$ and for all $x\in f_i(\oS_i)$ there is an open neighborhood $U_x$ of $x$ in $f_i(\oS_i)$ diffeomorphic to the open $i$-ball $B^i$ such that there is a diffeomorphism of stratifolds $\psi: r_i^{-1}(U_x)\isopfeil B^{i}\times F$, with $F$ a stratifold whose $0$-stratum $F_0$ is a single point. Here $r_i$ is the retract given by the collars.
 
A regular stratifold $S$ is called an \emph{Euler stratifold} if for all $x\in S$ these $F$ have the property that the complement of the $0$-stratum has even Euler characteristic: $\chi(F\setminus F_0)\equiv 0 \mod 2$.

An $n$-dimensional stratifold $T$ with boundary is called regular if its interior $S$ is regular, and Euler if its interior $S$ is Euler.
\end{Def}

We denote the homology theory defined by the class of Euler stratifolds by $Eh$. We obtain:

\begin{Prop} Let $H$ be a finite group. Then $$Eh^H_{n}(*)\cong \Zz$$ for all $n\geq 0$.
\end{Prop}
\begin{proof}
Let $[S]\in Eh^H_{n}(*)$. If $\chi(S)\equiv 0 \mod 2$, we can cone off $S$ as above, obtaining $[S]=[\partial CS]=0$~\cite[Lemma 7]{diplomarbeit}. If $\chi(S)\equiv 1 \mod 2$ we endow the one-point space $pt$ with the trivial $H$-operation, making it a $0$-dimensional $H$-stratifold, and then we cone off $S\cup pt$, which is possible since $\chi(S\cup pt)\equiv 0 \mod 2$. We obtain $[S\cup pt]=[\partial(C(S\cup pt))]=0$, and so $[S]=[pt]$. Since $\chi(\partial S)=0$ for all Euler stratifolds $S$, we know that $pt$ cannot be zero bordant~\cite[Satz 2]{diplomarbeit}. Thus $Eh^H_n(*)=\{0,[pt]\}\cong \Zz$. (Equivalently, we know that we have a surjection onto $Eh_{n}(*)\cong \Zz$.)   
\end{proof}

\begin{Cor} Let $G$ be a discrete group, and let $H$ be a finite subgroup of $G$. Then $$Eh^G_{n}(G/H)\cong \Zz$$ for all $n\geq 0$.
\end{Cor}
\begin{proof}
$Eh^G_{n}(G/H)\cong Eh^H_{n}(*) \cong \Zz$.
\end{proof}

Note that with this information, one can calculate $Eh^G_*(X)$ for all proper $G$-CW-complexes $X$ since these are built up of cells of the form $G/H\times D^i$ with $H$ a finite subgroup of $G$, $i\geq 0$. We can use the Mayer-Vietoris-sequence to compute $Eh^G_*(X)$ since we know $Eh^G_{n}(G/H)$ for all finite subgroups $H$ of $G$, for all $n\geq 0$. In particular, if $G$ is finite one knows $Eh^G_*(X)$ for all $G$-CW-complexes $X$.

\subsection{$\Cl=$ ``oriented singular homology stratifolds''}

An $n$-dimensional stratifold $S$ belongs to $\Cl$ if and only if its top-stratum $S^n$ is oriented, $G$ operates orientation preserving, and its ``codimension-1''-stratum $S^{n-1}$ is empty. The homology theory without $G$-action defined by these stratifolds is singular homology with $\Z$-coefficients, $H_*(-;\Z)$~\cite{kreck}. In the case of a $G$-action on the stratifolds, we obtain the following result:

\begin{Prop} Let $\Cl$ be the class of ``oriented singular homology stratifolds'' and let $H$ be a finite group. Then $$\H^H_{\Cl,n}(*)\cong \begin{cases} A(H) \text{ if }n=0\\ 0 \text{ else. }\end{cases}$$
Here $A(H)$ denotes the Burnside ring of $H$. 
\end{Prop}
\begin{proof}
Let $n>0$. For $[S]\in \H^H_{\Cl,n}(*)$, the cone $CS$ over $S$ again has an orientation on the top-stratum $S^n\times I$ and the codimension-1-stratum is $S^{n-1}\times I=\emptyset \times I = \emptyset$. We can define the action of $H$ on $CS=S\times [0,1]/S\times \{1\}$ as above by $h(x,t):=(hx,t).$ The $H$-stratifold $CS$ thus obtained gives a zero bordism of $S$: $[S]=[\partial CS]=0$. Thus for $n>0$ we have $\H^H_{\Cl,n}(*)=0$. 

Now let $n=0$ and let $[S]\in \H^H_{\Cl,0}(*)$. $S$ is an oriented $0$-dimensional stratifold, thus an oriented $0$-dimensional manifold, with $G$-action. The cone $CS$ cannot be used as a zero bordism, since it does not fulfill the requirement ``codimension-1-stratum empty'': its $0$-stratum is the cone point $\overline{S\times\{1\}}$. A bordism can only be given by an oriented $1$-dimensional stratifold with boundary whose $0$-stratum is empty - and this is nothing but a $1$-dimensional oriented smooth manifold with boundary. Thus we know that $\H^H_{\Cl,0}(*)=\Omega^H_0(*;\Z)$, and this is known to be $A(H)$~\cite{stong1970}.
\end{proof}

\begin{Cor} Let $G$ be a discrete group, and let $H$ be a finite subgroup of $G$. Then $$\H^G_{\Cl,n}(G/H)\cong \begin{cases} A(H) \text{ if }n=0\\ 0 \text{ else. }\end{cases}$$
Here $A(H)$ denotes the Burnside ring of $H$. 
\end{Cor}
\begin{proof}
$$\H^G_{\Cl,n}(G/H)\cong \H^H_{\Cl,n}(*)\cong \begin{cases} A(H) \text{ if }n=0\\ 0 \text{ else. }\end{cases}$$
\end{proof}

Note that with this information, one can calculate $\H^G_{\Cl,*}(X)$ for all proper $G$-CW-complexes $X$ since these are built up of cells of the form $G/H\times D^i$ with $H$ a finite subgroup of $G$, $i\geq 0$. We can use the Mayer-Vietoris-sequence to compute $\H^G_{\Cl,*}(X)$ since we know $\H^G_{\Cl,n}(G/H)$ for all finite subgroups $H$ of $G$, for all $n\geq 0$. In particular, if $G$ is finite one knows $\H^G_{\Cl,*}(X)$ for all $G$-CW-complexes $X$.

\subsection{$\Cl=$ ``non-oriented singular homology stratifolds''}

Now an $n$-dimensional stratifold $S$ belongs to $\Cl$ if and only if its ``codimension-1''-stratum $S^{n-1}$ is empty. There are no orientation requirements. Without $G$-action the homology theory defined by these stratifolds is singular homology with $\Zz$-coefficients, $H_*(-;\Zz)$~\cite{kreck}. In the case of a $G$-action on the stratifolds, we obtain:
\begin{Prop} Let $\Cl$ be the class of ``non-oriented singular homology stratifolds'' and let $H$ be a finite group. Then $$\H^H_{\Cl,n}(*)\cong \begin{cases} V \text{ if }n=0\\ 0 \text{ else, }\end{cases}$$
where $V$ is a $\Zz$-vector space with base $\{ H/K \}$, where $K$ belongs to a complete set of conjugacy class representatives of the collection of subgroups of $H$ having odd index in their normalizer.
\end{Prop}

\begin{proof}
Let $n>0$. As above, any stratifold $[S]\in \H^H_{\Cl,n}(*)$ can be coned off, since $CS^n=S^{n-1}\times I=\emptyset$, and we can extend the $G$-operation to $CS$. So $\H^H_{\Cl,n}(*)=0$. 

Let $n=0$ and let $ [S]\in \H^H_{\Cl,0}(*)$. Again a bordism cannot be given by the cone, only by a $1$-dimensional smooth manifold with boundary. Thus $\H^H_{\Cl,0}(*)=\N^H_0(*)$, and this is calculated in~\cite[Prop.~13.1]{stong1970} to be the $\Zz$-vector space $V$ described above.
\end{proof}

\Rem
In~\cite{stong1970} it is noted that there is a one to one correspondence between the collection of conjugacy classes of subgroups $K\leq H$ with $[N_H(K):K]$ odd and the collection of conjugacy classes of subgroups $L\leq H$ admitting no nontrivial homomorphism to $\Zz$. 

\begin{Cor} Let $G$ be a discrete group, and let $H$ be a finite subgroup of $G$. Then $$\H^G_{\Cl,n}(G/H)\cong \begin{cases} V \text{ if }n=0\\ 0 \text{ else, }\end{cases}$$
where $V$ is a $\Zz$-vector space with base $\{ H/K \}$, where $K$ belongs to a complete set of conjugacy class representatives of the collection of subgroups of $H$ having odd index in their normalizer.
\end{Cor}
\begin{proof}
$$\H^G_{\Cl,n}(G/H)\cong \H^H_{\Cl,n}(*)\cong \begin{cases} V \text{ if }n=0\\ 0 \text{ else. }\end{cases}$$
\end{proof}

Note that with this information, one can calculate $\H^G_{\Cl,*}(X)$ for all proper $G$-CW-complexes $X$ since these are built up of cells of the form $G/H\times D^i$ with $H$ a finite subgroup of $G$, $i\geq 0$. We can use the Mayer-Vietoris-sequence to compute $\H^G_{\Cl,*}(X)$ since we know $\H^G_{\Cl,n}(G/H)$ for all finite subgroups $H$ of $G$, for all $n\geq 0$. In particular, if $G$ is finite one knows $\H^G_{\Cl,*}(X)$ for all $G$-CW-complexes $X$.

\section{Possible further developments} \label{sec5}

One could continue this analysis for any other bordism class of stratifolds one is interested in. (The standard procedure would be to have ideas of how to calculate the bordism groups in the non-equivariant case, and to then generalize these to the equivariant setting.)

It would be more satisfying to solve the above problem completely also for discrete $G$ and non-proper $G$-CW-complexes $X$, but at the moment I do not see a promising approach. One might try to do something for actions of groups $G$ which are direct products of a finite group $F$ and an infinite group $Z$ without any nontrivial finite subgroups. Maybe one could apply the methods of Section~\ref{sec2} to the action of $Z$ (which is free since it is proper and $Z$ has no finite subgroups which could serve as isotropy groups other than the trivial group) and the methods of Section~\ref{sec4} to $F$ and then combine the results. All finitely generated abelian groups are of the type described above for which this method might work.

As already mentioned, it would also be interesting to consider restrictions on the $G$-operations on the stratifolds other than ``free and proper'' and ``proper''. (One has to be very careful in adapting the induction structure to these more general cases and probably needs additional requirements, somenthing like $\alpha(\F_H)\subseteq \F_G$ if one has the condition ``with isotropy in $\F$'', maybe more.)

It would also be interesting to extend the analysis to actions of higher dimensional Lie groups $G$. I think that the material in Sections~\ref{sec1} and~\ref{sec2} can be extended to this more general case, but this has to be carefully checked. One could define
\begin{eqnarray*}
\H^{G}_{\Cl,n}(X)= & \{ [S,f]\; | & S \text{ ($n + \dim G$)-dim.~cocomp.~$G$-stratifold in $\Cl$ } \\
& & \text{with free proper $G$-action, $f:S\to X$ $G$-equiv.}\}
\end{eqnarray*}
and see if all arguments used in the proofs are still valid. 

\bibliographystyle{alpha}

%\bibliography{equivstratis}

\end{document}